\documentclass[12pt]{article}
\usepackage{amscd,amssymb,amsmath}

\usepackage[dvips]{graphicx}
\usepackage[mathscr]{eucal}

\begin{document}
\newtheorem{proposition}{Proposition}
\newtheorem{lemma}{Lemma}
\newtheorem{corollary}{Corollary}
\newtheorem{theorem}{Theorem}

\centerline{\Large On a q-analog of a Sahi result\bf }
\centerline{\Large \bf}

\centerline {Olga Bershtein} \centerline{Institute for Low
Temperature Physics and Engineering, Kharkov, Ukraine}
\centerline{e-mail: bershtein@ilt.kharkov.ua}

\begin{abstract}
We obtain a $q$-analog of a well known Sahi result on the joint spectrum of
\hbox{$S(GL_n \times GL_n)$}-invariant differential operators with
polynomial coefficients on the vector space of complex $n \times
n$-matrices.

\medskip


\end{abstract}

\section{Introduction}

Start with recalling some well-known facts. Denote by
$\mathrm{Mat}_n$ the vector space of complex $n \times n$-matrices.
The group $K=S(GL_n \times GL_n)$ acts on $\mathrm{Mat}_n$ by
$$
(u,v)Z=uZv^{-1}, \quad (u,v) \in K, Z \in \mathrm{Mat}_n.
$$

This induces the natural $K$-actions in the vector spaces $\mathbb
C[\mathrm{Mat}_n]$ of holomorphic polynomials,
$\mathrm{D}[\mathrm{Mat}_n]$ of differential operators with constant
coefficients, and $\mathrm {PD}[\mathrm{Mat}_n]$ of differential
operators with polynomial coefficients. The well-known Hua theorem
claims that
$$
\mathbb C[\mathrm{Mat}_n]=\bigoplus_{\lambda \in \Lambda_n}
\mathbb{C} [\mathrm{Mat}_n]_{\bf{\lambda}},
$$
where $\Lambda_n=\{\lambda= (\lambda_1,\lambda_2,\ldots,\lambda_n)\in
\mathbb{Z}^n_+ | \lambda_1 \geq \lambda_2 \geq \ldots \geq
\lambda_n\}$ and $\mathbb{C}[\mathrm{Mat}_n]_{\bf{\lambda}}$ is a
simple finite dimensional $K$-module with the highest weight
\begin{equation}\label{h_weight}
(\lambda_1-\lambda_2,\ldots,\lambda_{n-1}-\lambda_n,2\lambda_n,
\lambda_{n-1}-\lambda_n,\ldots,\lambda_1-\lambda_2).
\end{equation}
Similarly, $\mathrm{D}[\mathrm{Mat}_n]=\bigoplus_{\lambda \in
\Lambda_n} \mathrm{D}_{\bf{\lambda}},$ with
$\mathrm{D}_{\bf{\lambda}} \cong
\mathbb{C}[\mathrm{Mat}_n]^*_{\bf{\lambda}}$ (we are using here the
standard pairing between polynomials and differential operators with
constant coefficients).

Let $y_{\nu}=\sum_i v_i w_i \in \mathrm {PD}[\mathrm{Mat}_n],$ where
$\{v_i\}$ is a basis in $\mathbb{C}[\mathrm{Mat}_n]_{\bf{\nu}}$, and
$\{w_i\}$ is the dual basis in $\mathrm{D}_{\bf{\nu}}$.
$y_{\nu}|_{\mathbb{C}[\mathrm{Mat}_n]_{\mathbf \lambda}}$ is a scalar
operator since $\mathbb{C}[\mathrm{Mat}_n]_{\mathbf \lambda}$ is a
simple $K$-module and $y_{\nu}$ is $K$-invariant. Sahi arranges
studying an explicit formula for these scalars \cite{
Sahi_to_Kostant}, \cite[Proposition 3.3]{KnopSahi}:
\begin{equation}\label{classic_y_spectr}
\mathrm{y}_{\nu}|_{\mathbb{C}[\mathrm{Mat}_n]_{\mathbf \lambda}}\quad
=\mathfrak{s}_{\nu}(\lambda_1 + n - 1,\lambda_2 + n -
2,\ldots,\lambda_{n-1}+1,\lambda_n),
\end{equation}
where the {\it factorial} Schur polynomial $\mathfrak{s}_{\nu}$
associated to a partition $\nu=(\nu_1,\ldots,\nu_n)$ is defined by
$$
\mathfrak{s}_\nu(x_1,x_2,\ldots,x_n) \quad=\quad
 \frac{\det\left(\prod\limits_{m=0}^{\nu_j+n-j-1}(x_i-m)
\right)_{1 \leq i,j \leq n}} {\prod\limits_{i<j} (x_i - x_j)},
$$
(see \cite{BL,Macdonald_Schur}).

This paper presents a $q$-analog of this formula.

\section{The main statement}\label{Pol_Mat}

Let $q \in (0,1)$. All algebras are assumed to be associative and
unital, and $\mathbb C$ is the ground field.

Recall that $U_q\mathfrak{sl}_{2n}$ is a Hopf algebra with generators
$\{E_i,\:F_i,\:K_i,\:K_i^{-1}\}_{i=1}^{2n-1}$ and the Drinfeld-Jimbo
relations \cite{Jantzen}

$$K_iK_j \,=\,K_jK_i,\quad K_iK_i^{-1}\,=\,K_i^{-1}K_i \,=\,1;$$
$$K_iE_i \,=\,q^{2}E_iK_i,\quad K_iF_i \,=\,q^{-2}F_iK_i;$$
$$K_iE_j \,=\,q^{-1}E_jK_i,\quad K_iF_j \,=\,qF_jK_i,\quad |i-j|=1;$$
$$K_iE_j \,=\,E_jK_i,\quad K_iF_j \,=\,F_jK_i,\quad |i-j| \,> \,1;$$
$$E_iF_j \,-\,F_jE_i \:=\:\delta_{ij}\frac{K_i-K_i^{-1}}{q-q^{-1}};$$
$$
E_i^2E_j\,-\,(q+q^{-1})E_iE_jE_i \,+\,E_jE_i^2 \:=\:0,\quad |i-j| \,=\,1;
$$
$$
F_i^2F_j \,-\,(q+q^{-1})F_iF_jF_i \,+\,F_jF_i^2\:=\:0, \quad |i-j|\,=\,1;
$$
$$E_iE_j-E_jE_i\,=\,F_iF_j-F_jF_i\,=\,0, \quad |i-j|\,\ne \,1.$$

The coproduct, the counit, and the antipode are defined as follows:
\begin{align*}
\triangle{E_j}&=E_j \otimes 1+K_j \otimes E_j,& \varepsilon(E_j)&=0,&
S(E_j)&=-K_j^{-1}E_j,\\ \triangle{F_j}&=F_j \otimes K_j^{-1}+1
\otimes F_j,& \varepsilon(F_j)&=0,& S(F_j)&=-F_jK_j,\\
\triangle{K_j}&=K_j \otimes K_j,& \varepsilon(K_j)&=1,&
S(K_j)&=K_j^{-1},\qquad j=1,\ldots,2n-1.
\end{align*}

Equip the Hopf algebra $U_q\mathfrak{sl}_{2n}$ with an involution $*$:
$$
(K_j^{\pm 1})^* = K_j^{\pm 1}, \quad E_j^* = \left\{
\begin{array}{rl}
   K_j F_j, & j \neq n, \\
  -K_j F_j, & j = n,
\end{array}
\right. \quad F_j^* = \left\{
\begin{array}{rl}
   E_j K_j^{-1}, & j \neq n, \\
  -E_j K_j^{-1}, & j = n.
\end{array}
\right.
$$

$U_q\mathfrak{su}_{n,n} \stackrel{\mathrm{def}}{=}
(U_q\mathfrak{sl}_{2n}, *)$ is a $*$-Hopf algebra. Denote by $U_q
\mathfrak{k} \subset U_q \mathfrak{sl}_{2n}$ the Hopf subalgebra
generated by $E_j, F_j, \, j \neq n,$ and $K_i, K_i^{-1},
i=1,\ldots,2n-1$.

Introduce a $*$-algebra $\mathrm{Pol}(\mathrm{Mat}_n)_q$, which is
one of the basic objects in the theory of quantum bounded symmetric
domains (see, for example, \cite{polmat}). First, introduce a
well-known quantum analog $\mathbb{C}[\mathrm{Mat}_n]_q$ of the
algebra $\mathbb C[\mathrm{Mat}_n]$ of holomorphic polynomials on the
matrix space (see, for example, \cite{Kl}, chap. 9.2). It is defined
by the generators $z_a^\alpha$, $a,\alpha=1,\ldots,n,$ and the
following relations
\begin{align}
& z_a^\alpha z_b^\beta-qz_b^\beta z_a^\alpha=0, & a=b \quad \& \quad
\alpha<\beta,& \quad \text{or}\quad a<b \quad \& \quad \alpha=\beta,
\label{zaa1}
\\ & z_a^\alpha z_b^\beta-z_b^\beta z_a^\alpha=0,& \alpha<\beta \quad
\&\quad a>b, & \label{zaa2}
\\ & z_a^\alpha z_b^\beta-z_b^\beta z_a^\alpha-(q-q^{-1})z_a^\beta
z_b^\alpha=0,& \alpha<\beta \quad \& \quad a<b.& \label{zaa3}
\end{align}

Similarly, denote by $\mathbb{C}[\overline{\mathrm{Mat}}_n]_q$ an
algebra with the generators $(z_a^\alpha)^*$, $a,\alpha=1,\dots,n$
and the defining relations
\begin{flalign}
& (z_b^\beta)^*(z_a^\alpha)^* -q(z_a^\alpha)^*(z_b^\beta)^*=0, & a=b
\; \& \; \alpha<\beta, & \quad \text{or} \quad a<b \; \& \;
\alpha=\beta, \label{zaa1*}
\\ & (z_b^\beta)^*(z_a^\alpha)^*-(z_a^\alpha)^*(z_b^\beta)^*=0, &
\alpha<\beta \; \& \; a>b, & \label{zaa2*}
\\ & (z_b^\beta)^*(z_a^\alpha)^*-(z_a^\alpha)^*(z_b^\beta)^*-
(q-q^{-1})(z_b^\alpha)^*(z_a^\beta)^*=0, & \alpha<\beta \; \& \; a<b.
& \label{zaa3*}
\end{flalign}

Let $\mathbb C[\mathrm{Mat}_n \oplus \overline{\mathrm{Mat}}_n]_q$ be
an algebra with the generators $z_a^\alpha$, $(z_a^\alpha)^*$,
$a,\alpha=1,\dots,n$, and the defining relations (\ref{zaa1}) --
(\ref{zaa3*}) and
\begin{equation}\label{z*z}
(z_b^\beta)^*z_a^\alpha = q^2 \sum\limits_{a',b'=1}^n
\sum\limits_{\alpha', \beta' = 1}^n R(b,a,b',a') R(\beta, \alpha,
\beta', \alpha') z_{a'}^{\alpha'} \left( z_{b'}^{\beta'} \right)^* +
(1-q^2) \delta_{ab} \delta^{\alpha \beta},
\end{equation}
where $\delta_{ab}$, $\delta^{\alpha \beta}$ are Kronecker symbols,
$$
R(j,i,j',i') = \left\{
\begin{array}{cl}
  q^{-1}, &\quad i \neq j\ \&\ j=j'\ \&\ i=i', \\
  1, &\quad i=j=i'=j', \\
  -(q^{-2}-1), &\quad i=j\ \&\ i'=j'\ \&\ i'>i, \\
  0, &\quad \mbox{otherwise}.
\end{array}
\right.
$$

Finally, let $\mathrm{Pol}(\mathrm{Mat}_n)_q$ denotes the $*$-algebra
$(\mathbb C[\mathrm{Mat}_n \oplus \overline{\mathrm{Mat}}_n]_q,*)$
with the involution given by \hbox{$*: z_a^\alpha \mapsto
(z_a^\alpha)^*$.}

It is very important for our purposes that
$\mathrm{Pol}(\mathrm{Mat}_n)_q$ is a $q$-analog of the algebra of
differential operators with polynomial coefficients
$\mathrm{PD}[\mathrm{Mat}_n]$ mentioned in the Introduction. Indeed,
the latter algebra is derivable from $\mathrm{Pol}(\mathrm{Mat}_n)_q$
via the change of generators $z_a^\alpha \rightarrow
(1-q^2)^{-1/2}z_a^\alpha$ and a subsequent formal passage to the
limit as $q \rightarrow 1$.

$\mathrm{Pol}[\mathrm{Mat}_n]_q$ can be equipped with a $U_q
\mathfrak{su}_{n,n}$-module algebra structure via such formulas (see
\cite[Sec. 9,10]{polmat}): for ${a, \alpha =1,\ldots,n}$
\begin{gather*}
K_n^{\pm 1}z_a^\alpha=
\begin{cases}
q^{\pm 2}z_a^\alpha,&a=n \;\&\;\alpha=n,
\\ q^{\pm 1}z_a^\alpha,&a=n \;\&\;\alpha \ne n \quad \mathrm{or}\quad a \ne
n \;\&\; \alpha=n,
\\ z_a^\alpha,&\mathrm{otherwise},
\end{cases}
\\ F_nz_a^\alpha=q^{1/2}\cdot
\begin{cases}
1,& a=n \;\& \;\alpha=n,
\\ 0,&\mathrm{otherwise},
\end{cases}\qquad
E_nz_a^\alpha=-q^{1/2}\cdot
\begin{cases}
q^{-1}z_a^nz_n^\alpha,&a \ne n \;\&\;\alpha \ne n,
\\ (z_n^n)^2,& a=n \;\&\;\alpha=n,
\\ z_n^nz_a^{\alpha},&\mathrm{otherwise},
\end{cases}
\end{gather*}
and for $k \ne n$
\begin{align*}
K_k^{\pm 1}z_a^\alpha&=
\begin{cases}
q^{\pm 1}z_a^\alpha,& k<n \;\&\;a=k \quad \mathrm{or}\quad k>n
\;\&\;\alpha=2n-k,
\\ q^{\mp 1}z_a^\alpha,& k<n \;\&\;a=k+1 \quad \mathrm{or}\quad k>n
\;\&\;\alpha=2n-k+1,
\\ z_a^\alpha,&\mathrm{ otherwise},
\end{cases}
\\ F_kz_a^\alpha&=q^{1/2}\cdot
\begin{cases}
z_{a+1}^\alpha,& k<n \;\&\;a=k,
\\ z_a^{\alpha+1},& k>n \;\&\;\alpha=2n-k,
\\ 0,&\mathrm{otherwise},
\end{cases},
\\ E_kz_a^\alpha&=q^{-1/2}\cdot
\begin{cases}
z_{a-1}^\alpha,& k<n \;\&\; a=k+1,
\\ z_a^{\alpha-1},& k>n \;\&\;\alpha=2n-k+1,
\\ 0,& \mathrm{otherwise}.
\end{cases}
\end{align*}

In the sequel we are using standard results on finite dimensional $U_q
\mathfrak{k}$-modules of type {\bf 1}, see \cite{Jantzen}. As a $U_q
\mathfrak{k}$-module, $\mathbb{C} [\mathrm{Mat}_n]_q = \bigoplus_{\lambda
\in \Lambda_n} \mathbb{C} [\mathrm{Mat}_n]_{q, \bf{\lambda}},$ with
$\mathbb{C} [\mathrm{Mat}_n]_{q, \bf{\lambda}}$ being a simple finite
dimensional $U_q \mathfrak{k}$-module with the highest weight given by
\eqref{h_weight}. Also,
\hbox{$\mathbb{C}[\overline{\mathrm{Mat}}_n]_q=\bigoplus_{\lambda \in
\Lambda_n} \mathbb{C} [\overline{\mathrm{Mat}_n}]_{q, \bf{\lambda}}$,} with
$\mathbb{C} [\overline{\mathrm{Mat}_n}]_{q, \bf{\lambda}} \approx
\mathbb{C} [\mathrm{Mat}_n]^*_{q, \bf{\lambda}}$ as $U_q
\mathfrak{k}$-modules. Then
$$
(\mathrm{Pol}(\mathrm{Mat}_n)_q)^{U_q \mathfrak{k}}=\bigoplus_{\nu
\in \Lambda_n}\mathrm{Pl}_{\nu}, \qquad \mathrm{Pl}_{\nu}=(\mathbb{C}
[\mathrm{Mat}_n]_{q, \bf{\nu}} \otimes \mathbb{C}
[\overline{\mathrm{Mat}_n}]_{q, \bf{\nu}})^{U_q \mathfrak{k}}, \qquad
\dim \mathrm{Pl}_{\nu}=1.
$$
Let $\{v_j\} \subset \mathbb{C} [\mathrm{Mat}_n]_{q, \bf{\nu}}$ be a
basis and $\{w_j\} \subset \mathbb{C} [\overline{\mathrm{Mat}_n}]_{q,
\bf{\nu}}$ the dual basis. Then $\sum_j v_jw_j \in
\mathrm{Pl}_{\nu}$. Introduce $q$-minors
$$
z_{\quad I}^{\wedge k\, J}\stackrel{\mathrm{def}}{=}\sum_{s \in
S_k}(-q)^{l(s)}z_{i_1}^{j_{s(1)}} z_{i_2}^{j_{s(2)}}\cdot \ldots
\cdot z_{i_k}^{j_{s(k)}},
$$
\centerline{$I=\{(i_1,i_2,\ldots,i_k)|1 \le i_1<i_2<\ldots<i_k \le
n\}$,} \centerline{$J=\{(j_1,j_2,\ldots,j_k)|1 \le
j_1<j_2<\ldots<j_k \le n\}$.} It can be verified easily that
$$
v_{\nu}=(z_{\quad \{1,\ldots,n\}}^{\wedge n\,
\{1,\ldots,n\}})^{\nu_n} \prod_{k=1}^{n-1}(z_{\quad
\{1,\ldots,k\}}^{\wedge k\, \{1,\ldots,k\}})^{\nu_k-\nu_{k+1}}
$$
is a highest weight vector of $\mathbb{C} [\mathrm{Mat}_n]_{q,
\bf{\nu}}$. Consider a basis $\{v_j\} \subset \mathbb{C}
[\mathrm{Mat}_n]_{q, \bf{\nu}}$ that contains $v_{\nu}$. The
isomorphism between $\mathbb{C} [\overline{\mathrm{Mat}_n}]_{q,
\bf{\nu}}$ and $\mathbb{C} [\mathrm{Mat}_n]_{q, \bf{\nu}}^*$
mentioned above can be chosen so that the dual basis $\{w_j\}$
contains $v_{\nu}^*$. Introduce $y_{\nu} \in \mathrm{Pl}_{\nu}$ by
$$y_{\nu}=\sum_j v_jw_j.$$

Denote by $\mathcal H$ a $\mathrm{Pol}(\mathrm{Mat}_n)_q$-module with
a generator $f_0$ and defining relations
$$(z_a^{\alpha})^*f_0=0, \qquad a, \alpha=1,\ldots,n,$$
and by $T_F$ the corresponding representation of
$\mathrm{Pol}(\mathrm{Mat}_n)_q$ in $\mathcal H.$ The statements of
the following proposition are proved in \cite{polmat}.
\begin{proposition}\label{fock}
\begin{enumerate}
\item $\mathcal {H}=\mathbb C[\mathrm{Mat}_n]_qf_0$.
 \item $\mathcal{H}$ is
a simple $\mathrm{Pol}(\mathrm{Mat}_n)_q$-module. \item There exists
a unique sesquilinear form $(\cdot,\cdot)$ on $\mathcal H$ with the
following properties:\\ i) $(f_0,f_0)=1$; ii) $(fv,w)=(v,f^*w)$ for
all $v,w \in \mathcal H$, $f \in \mathrm{Pol}(\mathrm{Mat}_n)_q$.
\item The form $(\cdot,\cdot)$ is positive definite. \item $T_F$ is a
faithful representation.
\end{enumerate}
\end{proposition}

So, $\mathcal{H}$ is a pre-Hilbert space, and $T_F$ is an irreducible
faithful $*$-representation.

Evidentally, $\mathcal{H}$ inherits the decomposition
\begin{equation}\label{decomp}
\mathcal{H}=\bigoplus_{\lambda \in \Lambda_n}
\mathcal{H}_{\bf{\lambda}}.
\end{equation}

\begin{proposition}{(D. Shklyarov)}
$y_{\nu}y_{\lambda}=y_{\lambda}y_{\nu}$ for all partitions $\nu \neq
\lambda$.
\end{proposition}
{\bf Proof.} Commutativity is deduced from the faithfulness of $T_F$
and the simplicity of the summands in \eqref{decomp}. \hfill
$\square$

As in the classical case, $T_F(y_{\nu})|_{\mathcal{H}_\lambda}$ are
scalar operators for all $\nu$ and $\lambda$. The main goal of this
paper is to obtain an explicit formula for the scalars
$T_F(y_{\nu})|_{\mathcal{H}_\lambda}$.

Recall the notation of the $q$-factorial Schur polynomials
\cite{KnopSahi}: for $\nu \in \Lambda_n$
$$
\mathfrak{s}_{\nu}(x_1,x_2,\ldots,x_n;q)=
\frac{\det\left(\prod\limits_{m=0}^{\nu_j+n-j-1}(x_i-q^m)) \right)_{1
\leq i,j \leq n}} {\prod\limits_{i<j} (x_i - x_j)}.
$$

\begin{theorem}\label{spectr_y_k}
For all partitions $\nu, \lambda \in \Lambda_n$
\begin{equation*}
T_F(y_{\nu})|_{\mathcal{H}_\lambda}= (-q)^{\sum_{i=1}^n \nu_i}
q^\mathrm{const}\,\mathfrak{s}_{\nu}
  (q^{2(\lambda_1 + n- 1)}, q^{2(\lambda_2 + n - 2)}, \ldots,
 q^{2(\lambda_{n-1} + 1)},  q^{2\lambda_n};q^2)
\end{equation*}
with $\mathrm{const}= -\sum_{i=1}^n \nu_i(\nu_i+2n-2i).$
\end{theorem}
This theorem is a natural generalization of the following result. Let
$\mathbf{1}^k\stackrel{\rm
def}{=}(\underbrace{1,\ldots,1}_k,0,\ldots,0)$.
\begin{theorem}\cite[Theorem 1]{BKV}\label{spectr_y_k}
For all $k=1,2,\ldots,n$ and all $\lambda \in \Lambda_n$
\begin{equation*}
T_F(y_{\mathbf{1}^k})|_{\mathcal{H}_\lambda}= (-q)^k
q^{-k^2-2k(n-k)}\,\mathfrak{s}_{\mathbf{1}^k}
  (q^{2(\lambda_1 + n- 1)}, q^{2(\lambda_2 + n - 2)}, \ldots,
 q^{2(\lambda_{n-1} + 1)},  q^{2\lambda_n};q^2).
\end{equation*}
\end{theorem}
First, we prove some auxiliary statements.
\begin{lemma}\label{gener}
The subalgebra $\mathrm{Pol}(\mathrm{Mat}_n)_q^{U_q \mathfrak{k}}$ is
generated by the elements $y_{\mathbf 1^k}$.
\end{lemma}
{\bf Proof.} Equip the $U_q \mathfrak{k}$-module
$(\mathrm{Pol}(\mathrm{Mat}_n)_q)^{U_q \mathfrak{k}}$ with the
natural grading
$$
(\mathrm{Pol}(\mathrm{Mat}_n)_q)^{U_q
\mathfrak{k}}=\bigoplus_{j=0}^\infty
(\mathrm{Pol}(\mathrm{Mat}_n)_q)^{U_q \mathfrak{k}}_j, \qquad
(\mathrm{Pol}(\mathrm{Mat}_n)_q)^{U_q \mathfrak{k}}_j=\bigoplus_{\nu
\in \Lambda_n, |\nu|=j} \mathrm{Pl}_{\nu}.
$$
Hence $\dim (\mathrm{Pol}(\mathrm{Mat}_n)_q)^{U_q
\mathfrak{k}}_j=\#\{\nu \in \Lambda_n, |\nu|=j\}$.

It follows from Theorem \ref{spectr_y_k} that monomials
$\{y_{\mathbf 1^1}^{a_1}y_{\mathbf 1^2}^{a_2}...y_{\mathbf
1^n}^{a_n}\}$ are linear independent. This fact allows one to denote
by $\mathbb C[y_{\mathbf 1^1},y_{\mathbf 1^2},...,y_{\mathbf 1^n}]
\subset (\mathrm{Pol}(\mathrm{Mat}_n)_q)^{U_q \mathfrak{k}}$ the
subalgebra in generated by the elements $y_{\mathbf 1^1},y_{\mathbf
1^2},...,y_{\mathbf 1^n}$. It is easy to see that $\deg y_{\mathbf
1^k}=k$. Hence $\mathbb C[y_{\mathbf 1^1},y_{\mathbf
1^2},...,y_{\mathbf 1^n}] \bigcap
(\mathrm{Pol}(\mathrm{Mat}_n)_q)^{U_q \mathfrak{k}}_j$ is a linear
span of $\{y_{\mathbf 1^1}^{a_1}y_{\mathbf 1^2}^{a_2}...y_{\mathbf
1^n}^{a_n}|a_1+2a_2+...+na_n=j\}$, and
$$
\dim(\mathbb C[y_{\mathbf 1^1},y_{\mathbf 1^2},...,y_{\mathbf 1^n}]
\bigcap (\mathrm{Pol}(\mathrm{Mat}_n)_q)^{U_q
\mathfrak{k}}_j)=\#\{a_1,...,a_n \in \mathbb
Z_+|a_1+2a_2+...+na_n=j\}.
$$
So, $\dim(\mathbb C[y_{\mathbf 1^1},y_{\mathbf 1^2},...,y_{\mathbf
1^n}] \bigcap (\mathrm{Pol}(\mathrm{Mat}_n)_q)^{U_q
\mathfrak{k}}_j)=\dim (\mathrm{Pol}(\mathrm{Mat}_n)_q)^{U_q
\mathfrak{k}}_j,$ and $(\mathrm{Pol}(\mathrm{Mat}_n)_q)^{U_q
\mathfrak{k}}$ is generated by $y_{\mathbf 1^1},y_{\mathbf
1^2},...,y_{\mathbf 1^n}$. \hfill $\square$

The next statements concern symmetric polynomials and vanishing
conditions in the spirit of papers \cite{Knop-Capelli, Sahi2} (for
their classical analogs, see \cite{KnopSahi,KostSahi1,KostSahi2}).
Recall some notations from \cite{Knop-Capelli}. Fix non-zero real
numbers $q$ and $t$. For every $\lambda \in \Lambda_n$ we define
$\bar\lambda=(q^{\lambda_1},
q^{\lambda_2}t^{-1},...,q^{\lambda_n}t^{-n+1})$. We use the following
short notation: $|\lambda|=\sum_{i=1}^n \lambda_i$ for $\lambda \in
\Lambda_n.$ Also, let $m_\lambda$ be the monomial symmetric
polynomial that corresponds to $\lambda \in \Lambda_n$. Recall that
the usual order on $\Lambda_n$: for $\lambda, \mu \in \Lambda_n$ we
say $\lambda \geq \mu$ if $\lambda_1+...+\lambda_i \geq
\mu_1+...+\mu_i$ for all $i=1,...,n$.

\begin{proposition}\cite[Theorem 2.4]{Knop-Capelli}\label{Knop}
For every $\lambda \in \Lambda_n$ there exists a unique symmetric
polynomial $P_{\lambda}(z;q,t)$ in $n$ variables such that
$P_{\lambda}(\bar\mu;q,t)=0$ for all $\mu \in \Lambda_n, |\mu| \leq
|\lambda|, \mu \neq \lambda$, and which has an expansion
$P_{\lambda}(z;q,t)=\sum_{\mu \leq \lambda}p_{\lambda\mu}m_\mu(z)$
with $p_{\lambda\lambda}=1$.
\end{proposition}
\begin{proposition}\cite[Proposition 2.8]{Knop-Capelli}\label{Knop1}
$P_\lambda(z;q,q)=q^{-(n-1)|\lambda|}\mathfrak{s}_\lambda(q^{n-1}z;q)$.
\end{proposition}
Now we can prove
\begin{lemma}\label{l_2}
For any partition $\nu$ there exists a constant $c_{\nu}$ such that
for all $\lambda \in \Lambda_n$
\begin{equation*}
T_F(y_{\nu})|_{\mathcal{H}_\lambda}= c_{\nu}\mathfrak{s}_{\nu}
  (q^{2(\lambda_1 + n- 1)}, q^{2(\lambda_2 + n - 2)}, \ldots,
 q^{2(\lambda_{n-1} + 1)},  q^{2\lambda_n};q^2).
\end{equation*}
\end{lemma}
{\bf Proof.} It follows from Theorem \ref{spectr_y_k} and Lemma
\ref{gener}, that for an arbitrary partition $\nu$ there exists a
symmetric polynomial of degree $|\nu|$ in $n$ variables
$x_1,...,x_n$, such that the eigenvalues
$T_F(y_{\nu})|_{\mathcal{H}_{\lambda}}$ are just the values of the
polynomial at $x_1=q^{2(\lambda_1+n-1)},\ldots,x_n=q^{2\lambda_n}$.

Let $\delta=(n-1,...,1,0)$, $q^{2(\mu + \delta)}\stackrel{\rm
def}{=}(q^{2(\mu_1 + n- 1)}, q^{2(\mu_2 + n - 2)}, \ldots,
q^{2(\mu_{n-1} + 1)}, q^{2\mu_n})$ for any $\mu \in \Lambda_n$.
Propositions \ref{Knop} and \ref{Knop1} claim that
$\mathfrak{s}_{\nu}(x_1,...,x_n;q^2)$ is a unique (up to a constant
multiplier) symmetric polynomial of degree $|\nu|$ with
$\mathfrak{s}_{\nu}(q^{2(\mu+\delta)};q^2)=0$ for all $\mu \in
\Lambda_n, |\mu| \leq |\nu|, \mu \neq \nu$.

To finish the proof, one should investigate zeros of $T_F(y_{\nu})$
to conclude the proof (cf. the proof in Sahi's paper
\cite{Sahi_to_Kostant}). We claim that
$$
T_F(y_{\nu})|_{\mathcal{H}_{\lambda}}=0 \quad \text{for $|\nu| \leq
|\lambda|$} \quad \text{unless} \quad \nu = \lambda.
$$ Indeed, it suffices to prove that $T_F(y_{\nu})(v_{\lambda}f_0)=0$
for partitions $\nu$ and $\lambda$ such that $|\nu| \leq |\lambda|$,
$\nu \neq \lambda$. Recall that $y_{\nu}=\sum v_jw_j$, where $\{v_j\}
\subset \mathbb{C} [\mathrm{Mat}_n]_{q, \bf{\nu}}$ contains $v_{\nu}$
and $\{w_j\} \subset \mathbb{C} [\overline{\mathrm{Mat}_n}]_{q,
\bf{\nu}}$ contains $v_{\nu}^*$. It follows from the commutation
relations \eqref{z*z} that $T_F(w_j)(v_{\lambda}f_0)=0$ unless
$w_j=v_{\lambda}^*$. So, $T_F(y_{\nu})(v_{\lambda}f_0)=0$ unless
$\nu=\lambda$. \hfill $\square$

Introduce the notations: $\lambda-a
\mathbf{1}^n=(\lambda_1-a,...,\lambda_n-a)$. The next proposition
completes the proof of Theorem 1.
\begin{proposition}
$c_{\nu}=(-q)^{|\nu|} q^{-\sum \limits_{i=1}^n \nu_i(\nu_i+2n-2i)}.$
\end{proposition}
{\bf Proof.} Let us compare $T_F(y_{\nu})(v_{\nu}f_0)$ and
$\mathfrak{s}_{\nu}(q^{2(\nu + \delta)};q^2)$. By Proposition
\ref{4},
$$
T_F(y_{\nu})(v_{\nu}f_0)=(\mathrm {det}_q \mathbf z)^{\nu_n}
T_F(y_{\nu-\nu_n\mathbf{1}^n})T_F((\mathrm{det}_q \mathbf
z)^*)^{\nu_n})v_{\nu}f_0.
$$
By Theorem \ref{spectr_y_k},
$$
T_F((\mathrm{det}_q \mathbf
z)^*)^{\nu_n})v_{\nu}f_0=(-1)^{n\nu_n}q^{-n(n-1)\nu_n}\prod_{i=0}^{\nu_n-1}
\mathfrak{s}_{\mathbf{1}^n}(q^{2(\nu+\delta-i\mathbf{1}^n)};q^2)
v_{\nu-\nu_n\mathbf{1}^n}f_0.
$$
We proceed by induction in $n$. For $n=1$ the statement follows from
the last identity and Lemma \ref{schur_1} from the next section.

Let $n>1$. One can rewrite some of the commutation relations
\eqref{z*z} more explicitly:
\begin{align*}
& (z_n^\alpha)^*z_b^\beta = q \sum\limits_{\alpha', \beta' = 1}^n
R(\alpha, \beta, \alpha', \beta') z_b^{\beta'} (z_n^{\alpha'})^*
\qquad & \text{for} \qquad b<n,
\\& (z_a^n)^*z_b^\beta =q \sum\limits_{a',b'=1}^n R(a,b,a',b')
z_{b'}^{\beta} (z_{a'}^n)^* \qquad & \text{for} \qquad \beta<n.
\end{align*}
Hence, $T_F((z_a^\alpha)^*)v_{\nu-\nu_n\mathbf{1}^n}f_0=0$ for $a=n$
or $\alpha=n$. Denote by $T'_F$ the faithful representation of
$\mathrm{Pol}(\mathrm{Mat}_{n-1})_q$ in the vector space
$\mathcal{H'}$ defined by a single generator $f_0'$ and the relations
$(z_a^{\alpha})^*f_0'=0$, for $a, \alpha=1,...,n-1$. Thus,
$$
T_F(y_{\nu-\nu_n\mathbf{1}^n})|_{\mathcal{H}_{\nu-\nu_n\mathbf{1}^n}}=
T'_F(y_{\tau})|_{\mathcal{H'}_{\tau}},
$$
where $\tau=(\nu_1-\nu_n,...,\nu_{n-1}-\nu_n)$. Let
$\delta'=(n-2,\ldots,2,1,0)$. Hence, by the inductive assumption,
\begin{align*}
T_F(y_{\nu-\nu_n\mathbf{1}^n})v_{\nu-\nu_n\mathbf{1}^n}f_0=
(-q)^{|\tau|} q^{-\sum_{i=1}^{n-1}\tau_i(\tau_i+2n-2-2i)}
\mathfrak{s}_{\tau}(q^{2(\tau+\delta')};q^2)v_{\nu-\nu_n\mathbf{1}^n}f_0.
\end{align*} Now the required statement follows from Lemmas
\ref{schur_1}, \ref{schur_2} of the next section and the following
computation
\begin{align*}
&T_F(y_{\nu})v_{\nu}f_0=
(-1)^{n\nu_n}q^{-n(n-1)\nu_n}\prod_{i=0}^{\nu_n-1}
\mathfrak{s}_{\mathbf{1}^n}(q^{2(\nu+\delta-i\mathbf{1}^n)};q^2)
(\mathrm {det}_q \mathbf z)^{\nu_n}
T_F(y_{\nu-\nu_n\mathbf{1}^n})v_{\nu-\nu_n\mathbf{1}^n}f_0=
\\ &(-1)^{n\nu_n}q^{-n(n-1)\nu_n}\prod_{i=0}^{\nu_n-1}
\mathfrak{s}_{\mathbf{1}^n}(q^{2(\nu+\delta-i\mathbf{1}^n)};q^2)
(-q)^{|\tau|} q^{-\sum_{i=1}^{n-1}\tau_i(\tau_i+2n-2-2i)}
\mathfrak{s}_{\tau}(q^{2(\tau+\delta')};q^2)v_{\nu}f_0=
\\& (-q)^{|\nu|} q^{-\sum_{i=1}^n\nu_i(\nu_i+2n-2i)}
\mathfrak{s}_{\nu}(q^{2(\nu+\delta)};q^2) v_{\nu}f_0. \hfill \square
\end{align*}

\begin{proposition}\label{4}
$y_{\nu}=(\det_q \mathbf z)^{\nu_n} y_{\nu-\nu_n\mathbf{1}^n}((\det_q
\mathbf z)^*)^{\nu_n}$.
\end{proposition}
{\bf Proof.} It is obvious that $(\det_q \mathbf z)^{\nu_n}
y_{\nu-\nu_n\mathbf{1}^n}((\det_q \mathbf z)^*)^{\nu_n} \in
\mathrm{Pl}_{\nu}=\mathbb C \cdot y_{\nu}$, the statement follows
from an explicit computation of the coefficient of
$v_{\nu}v_{\nu}^*$.\hfill $\square$

\section {q-factorial Schur functions}

This section contains auxiliary statements which we used above. As
usual, \hbox{$(a)_n \stackrel{\rm def}{=} \prod
\limits_{i=0}^{n-1}(a-q^{2i})$.}

\begin{lemma}\label{schur_1}
For any partition $\nu \in \Lambda_n$, such that $\nu_n
>0$, one has $$\mathfrak{s}_{\nu}(q^{2\nu+2\delta};q^2)= q^{2|\nu|-2n}
\mathfrak{s}_{\mathbf{1}^n}(q^{2(\nu+\delta)};q^2)
\mathfrak{s}_{\nu-\mathbf{1}^n} (q^{2(\nu+\delta-\mathbf{1}^n)}
;q^2).$$
\end{lemma}
{\bf Proof.} The proof reduces to the explicit computation
$\mathfrak{s}_{\nu}(q^{2\nu+2\delta};q^2)=$
\begin{align*}
&=\prod_{1 \leq i \leq j \leq n}
\frac{1}{q^{2\nu_i+2n-2i}-q^{2\nu_j+2n-2j}}
\begin{vmatrix}
(q^{2\nu_1+2n-2})_{\nu_1+n-1} & 0 & ... & 0
\\ (q^{2\nu_1+2n-2})_{\nu_2+n-2} &
(q^{2\nu_2+2n-4})_{\nu_2+n-2} & ... & 0
\\ ...& ...& ...& ...
\\ (q^{2\nu_1+2n-2})_{\nu_n} &
(q^{2\nu_2+2n-4})_{\nu_n} & ... & (q^{2\nu_n})_{\nu_n}
\end{vmatrix}=
\\& = \prod_{i=1}^n(q^{2\nu_i+2n-2i}-1)q^{2\nu_i+2n-4i-2} \prod_{1
\leq i \leq j \leq n} \frac{1}{q^{2\nu_i+2n-2i-2}-q^{2\nu_j+2n-2j-2}}
\cdot
\\& \cdot \begin{vmatrix} (q^{2\nu_1+2n-4})_{\nu_1+n-2} & 0 &
... & 0
\\ (q^{2\nu_1+2n-4})_{\nu_2+n-3} &
(q^{2\nu_2+2n-6})_{\nu_2+n-3} & ... & 0
\\ ...& ...& ...& ...
\\ (q^{2\nu_1+2n-4})_{\nu_n-1} &
(q^{2\nu_2+2n-6})_{\nu_n-1} & ... & (q^{2\nu_n-2})_{\nu_n-1}
\end{vmatrix}
\\&
 = q^{2|\nu|}q^{-2n}\mathfrak{s}_{\mathbf{1}^n}(q^{2\nu+2\delta};q^2)
\mathfrak{s}_{\nu-\mathbf 1^n}(q^{2(\nu+\delta-\mathbf{1}^n)};q^2).
\hfill \square
\end{align*}

\begin{lemma}\label{schur_2}
For a partition $\nu \in \Lambda_n$ with $\nu_n=0$ one has
$\mathfrak{s}_{\nu}(q^{2(\nu+\delta)};q^2)=q^{2|\nu|}
\mathfrak{s}_{\nu'}(q^{2(\nu'+\delta')};q^2),$ where
$\nu'=(\nu_1,\ldots,\nu_{n-1})$.
\end{lemma}
{\bf Proof.} The proof is managed by the explicit computation
\begin{align*}
& \mathfrak{s}_{\nu}(q^{2(\nu+\delta)};q^2)= \prod_{1 \leq i \leq j
\leq n} (q^{2\nu_i+2n-2i}-q^{2\nu_j+2n-2j})^{-1}
\begin{vmatrix}
(q^{2\nu_1+2n-2})_{\nu_1+n-1} & 0 & ... & 0
\\ (q^{2\nu_1+2n-2})_{\nu_2+n-2} &
(q^{2\nu_2+2n-4})_{\nu_2+n-2} & ... & 0
\\ ...& ...& ...& ...
\\ 1 &
1 & ... & 1
\end{vmatrix}
\\ & = \prod_{1 \leq i \leq j < n}
(q^{2\nu_i+2n-2i}-q^{2\nu_j+2n-2j})^{-1}
\begin{vmatrix}
\frac{(q^{2\nu_1+2n-2})_{\nu_1+n-1}}{q^{2\nu_1+2n-2}-1} & 0 & ... & 0
\\ \frac{(q^{2\nu_1+2n-2})_{\nu_2+n-2}}{q^{2\nu_1+2n-2}-1} &
\frac{(q^{2\nu_2+2n-4})_{\nu_2+n-2}}{q^{2\nu_2+2n-4}-1} & ... & 0
\\ ...& ...& ...& ...
\\ \frac{(q^{2\nu_1+2n-2})_{\nu_{n-1}}}{q^{2\nu_1+2n-2}-1} &
\frac{(q^{2\nu_2+2n-4})_{\nu_{n-1}}}{q^{2\nu_2+2n-4}-1} & ... &
\frac{(q^{2\nu_{n-1}+2})_{\nu_{n-1}}}{q^{2\nu_{n-1}+2}-1}
\end{vmatrix}
\\ & = q^{2|\nu|}\mathfrak{s}_{\nu'}(q^{2\nu'+2\delta'};q^2).  \hfill \square
\end{align*}

\section{Acknowledgements}

The author thanks to L. Vaksman for helping with the proof of Lemma
\ref{l_2} and constant attention to her work. Also thanks are due to
D.Shklyarov. At last, the author would like to express her gratitude
to a referee for many useful remarks.


\begin{thebibliography}{15}

\bibitem{BKV} O. Bershtein, Ye. Kolisnyk, L. Vaksman, {\it On a
q-analog of the Wallach-Okounkov formula}, -- Lett. in Math.Phys.
 78 No. 1 (2006), pp.97-109.

\bibitem{BL} L.Biedenharn, J.Louck, {\it A new class of symmetric
polynomials defined in terms of tableaux}, -- Adv. in Appl. Math. 10
(1989), pp.396-438.

\bibitem{Jantzen} J.C. Jantzen, {\it Lectures on Quantum Groups} --
Amer. Math. Soc., Providence RI, (1996).

\bibitem{Kl} A. Klimyk, K. Schm\"udgen, {\it Quantum Groups and Their
Representations} -- Springer, Berlin, (1997).

\bibitem{KnopSahi} F. Knop, S. Sahi, {\it Difference equations and
symmetric polynomials defined by their zeros}, -- IMRN 10 (2000),
pp.473-486.

\bibitem{Knop-Capelli} F. Knop, {\it Symmetric and non-symmetric
quantum Capelli polynomials}, -- Comment. Math. Helv. 72 (1997),
pp.84-100.

\bibitem{KostSahi1} B.Kostant, S.Sahi, {\it The Capelli identity,
tube domains and the generalized Laplace transform}, -- Adv.Math. 87
(1991), pp.71-92.

\bibitem{KostSahi2} B.Kostant, S.Sahi, {\it Jordan algebras and
Capelli identities}, -- Inv.Math. 112 (1993), pp.657-664.

\bibitem{Macdonald_Schur} I. Macdonald, {\it Schur functions: theme
and variations}, -- I.R.M.A. Strasbourg (1992) 498/S-27, pp.5-39.

\bibitem{Sahi2} S.Sahi, {\it Interpolation, integrality and a
generalization of Macdonald's polynomials}, -- Int.Math.Res.Not. 10
(1996), pp.457-471.

\bibitem{Sahi_to_Kostant} S. Sahi, {\it The Spectrum of certain
invariant differential operators associated to a Hermitian symmetric
space}, -- in Lie Theory and Geometry, ed. J.-L.Brylinsky,
R.Brylinsky, V.Guillemin, V.Kac (1994), pp. 569-576.

\bibitem{qabsd} D. Shklyarov, S. Sinel'shchikov, L. Vaksman, {\it
q-Analogs of some bounded symmetric domains}, -- Czech. J. of Phys.
50 No.1 (2000), pp. 175-180.

\bibitem{polmat} D. Shklyarov, S. Sinel'shchikov, L. Vaksman, {\it
Fock representations and quantum matrices}, -- International J.Math
15 No.9 (2000), pp.1-40.
\end{thebibliography}
\end{document}